\documentclass[a4paper,12pt]{article} 

\usepackage{amsmath,amssymb}
\usepackage[usenames]{color}
\usepackage{graphicx}
\usepackage{amsfonts}
\usepackage{comment}
\usepackage{pxrubrica}

\newtheorem{thm}{Theorem}
\newtheorem{prop}{Proposition}
\newtheorem{lmm}{Lemma}
\newtheorem{defn}{Definition}

\numberwithin{equation}{section}

\newcommand{\qed}{\hfill\ensuremath{\square}\\}%

\newcommand{\pr}{\par \vspace{0mm} \noindent {\bf [Proof]} \quad}
\newcommand{\prend}{\hfill \qed}

\newcommand{\1}{{\bf 1}} 

\newcommand{\wt}{{\rm wt}}

\newcommand{\End}{{\rm End}}

\newcommand{\Com}{{\rm Com}}

\newcommand{\CF}{{\cal F}}

\newcommand{\bC}{{\mathbb C}} 
 
\newcommand{\bZ}{{\mathbb Z}} 
\newcommand{\bN}{{\mathbb N}}

\newcommand{\bR}{{\mathbb R}}

\definecolor{skyblue}{rgb}{0.5,0.5,1}
\definecolor{fadegreen}{rgb}{0.9,1,0.9}

\begin{document}
\title{Rationality of holomorphic vertex operator algebras}
\author{
\begin{tabular}{c}
Masahiko Miyamoto
\footnote{Partially supported
by the Grants-in-Aids
for Scientific Research, No.21K03195 and 18K18708, The Ministry of Education,
Science and Culture in Japan
and Academia Sinica in Taiwan} \cr
Institute of Mathematics, University of Tsukuba
\end{tabular}}
\date{}
\maketitle

\begin{abstract}
We prove that if $V$ is a unitary simple holomorphic vertex operator algebra of CFT-type, then $V$ is rational, that is, all $\bN$-gradable
$V$-modules are direct sums of copies of $V$. 
\end{abstract}

\section{Introduction}
A concept of vertex operator algebras (shortly VOA) was introduced by Borcherds \cite{B} to solve the moonshine conjecture in finite group theory.
It is now considered to be a rigorous algebraic definition of a chiral algebra
in conformal field theory.
In the last two decades, there has been quite good progress on the classification of holomorphic VOAs $V$ of CFT-type and with central charge $24$ (for example, see \cite{EMS} and
\cite{LM}).
In these works, to use modular invariance properties and to get the fact that weight one space is a reductive Lie algebra, most papers were assuming strongly regular conditions on $V$, that is,
positivity (or unitarity) of an invariant bilinear form, rationality (i.e., semisimplicity for $\bN$-gradable modules), and $C_2$-cofiniteness on $V$ (or finiteness of simple (weak) modules).
For the definition of unitarity, see \S 2.
Compared with the classification of unimodular positive definite even lattices \cite{N}, it is natural to
assume positivity over $\bR$ (or unitarity over $\bC$) for an invariant bilinear form on $V$, but is it necessary to
assume the rationality and $C_2$-cofiniteness?
For some holomorphic VOAs of CFT-type with some stronger assumptions,
their rationalities are already known, see \cite{DLM}.
In this paper, we will show that rationality comes from only the positivity (or unitary) of an invariant bilinear form of $V$ for a holomorphic VOA $V$ of CFT-type without any additional
assumptions. We note that every holomorphic VOA has a
symmetric invariant bilinear form $\langle\, ,\,\rangle$ by \cite{L}.

\begin{thm}\label{Main}
If $V=\coprod_{n=0}^{\infty}V_n$ is a unitary simple holomorphic VOA of CFT-type, then $V$ is rational, that is,
all $\bN$-gradable $V$-modules are direct sums of copies of $V$,
where holomorphic means that $V$ itself is the only one
simple $\bN$-gradable $V$-module.
\end{thm}

The author expresses great thanks to T.~Abe, C.~H.~Lam, and
H.~Yamauchi for suggesting many improvements to this note and also checking the proof. He thanks the staff of Academia Sinica in Taiwan for their hospitality. The author also expresses thanks to the referee for his/her advice.

\section{Two propositions and preliminary results}
Throughout this paper, $(V, Y, \1, \omega)$ denotes a
simple holomorphic VOA over $\bC$. 
In this paper, a $V$-module $M$ is called $\bN$-gradable if there is an $\bN$-grading structure $M=\oplus_{m=0}^{\infty}M_{(m)}$ such that
$v_m(M_{(n)})\subseteq M_{(n+\wt(v)-m-1)}$ for all $v\in V_{\wt(v)}$. We call $M_{(0)}$ the top level of $M=\oplus_{m=0}^{\infty}M_{(m)}$.
We also use $\wt(u)$ to denote an eigenvalue of $L(0)$ on an eigenvector $u\in M$ and
$M_t$ denotes the eigenspace of $L(0)$ with an eigenvalue $t$ for a $V$-module $M$.
For a simple $\bN$-gradable module $M$, it is easy to see that there is $-t\in \bC$ such that $M_{(n)}=M_{n-t}$ for all $n$.
Let $A(V)$ be the Zhu algebra of a VOA $V$ introduced by Zhu \cite{Z}. Then for an $\bN$-gradable $V$-module
$M=\oplus_{m=0}^{\infty}M_{(m)}$, $M_{(0)}$ is an $A(V)$-module.
One of the crucial results of Zhu's theory we will use is that every simple $A(V)$-module is a top level of an $\bN$-gradable simple $V$-module. 
Furthermore, for any $A(V)$-module $Z$ and a $V$-module $M$ generated from $Z$, we can define a canonical $\bN$-grading
$M=\coprod_{m=0}^{\infty}M_{(m)}$ on $M$ such that the top level $M_{(0)}$ of $M$ is $Z$. In particular, since $A(V)$ is a an $A(V)$-module, if $V$ is rational, then 
$A(V)$ is a direct sum of simple $A(V)$-modules. Namely, $A(V)$ is semisimple. 

If $V$ is simple and holomorphic, then $V=\coprod_{n=-t}^{\infty}V_n$ (with $\dim V_{-t}\not=0$) is the only simple $\bN$-gradable module 
and $A(V)/J(A(V))=\End(V_{-t})\cong {\rm Mat}_{s,s}(\bC)$, where $s=\dim V_{-t}$ and
$J(A(V))$ denotes the Jacobson radical of $A(V)$. 
We note that if $A(V)$ is Artinian, then most results in this section are well-known, but we have to treat the other cases, for example, 
the case where $\dim J(A(V))=\infty$.
We will find that $J(A(V))=0$ implies the rationality of $V$ (see Lemma 2).
Therefore, we study $J(A(V))$ for a while.
Since $A(V)/J(A(V))\cong \End(V_{-t})$, $J(A(V))$ is the unique maximal both-side ideal and is the intersection
of all maximal left ideals and coincides with the intersection of all maximal right ideals. 
In particular, for $\alpha\in J(A(V))$, $1-\alpha$ is a unit element of $A(V)$.
As an application, we have the following elementary lemma:

\begin{lmm}
For an $\bN$-graded $V$-module $M=\coprod_{m=0}^{\infty}M_{(m)}$ and any $0\not=a\in M_{(0)}$, we have $J(A(V))a \not=A(V)a$.
\end{lmm}

\pr Suppose false, then there is $\alpha\in J(A(V))$ such that $\alpha a=a$,
that is, $(1-\alpha)a=0$, which contradicts the fact that $1-\alpha$ is a unit element. 
\prend

We recall the notion of unitary VOAs from \cite{DLin14}.

\begin{defn}[\cite{DLin14}]
Let $(V, Y, \1, \omega)$ be a VOA over $\bC$ and
let $\phi:V\to V$ be an anti-linear involution of $V$ $($i.e, $\phi(\lambda v)=\bar{\lambda}\phi(u)$, $\phi(\1) = \1$, $\phi(\omega) =\omega$,
$\phi(u_nv) = \phi(u)_n\phi(v)$ for any $u, v \in V$, $n\in \bZ$, and
$\phi$ has order $2$ $)$.
Then $(V,\phi)$ is said to
be unitary if there exists a positive-definite Hermitian form $(\cdot,\cdot)_V:V\times V\to \bC$, which
is $\bC$-linear on the first vector and anti-$\bC$-linear on the second vector, such that the following
invariant property holds of any $a,u,v\in V$:
$$(Y (e^{zL(1)}(-z^{-2})^{L(0)}a, z^{-1})u, v)_V=(u, Y(\phi(a), z)v)_V, \eqno{(2.1)}$$
where $L(n)$ is defined by $Y(\omega,z)=\sum_{n\in \bZ} L(n)z^{-n-2}$.
\end{defn}

In our case, since $V$ is a holomorphic VOA, $V$ has a non-singular invariant bilinear form $\langle \cdot,\cdot\rangle$ satisfying
(2.1) by \cite{L}. Hence if we put $(u,v)_V=\langle u,\phi(v)\rangle$ for $v,u\in V$, then $(\cdot,\cdot)_V$ becomes a desired Hermitian form.
Therefore, $V$ is unitary if and only
if there is an anti-linear involution $\phi$ of $V$ such that $\langle v,\phi(v)\rangle>0$ for all $0\not=v\in V$.
We fix an anti-linear involution $\phi$ of $V$.

\begin{defn}[\cite{DLM}]
Let $M$ be a (weak) $V$-module.
We call an element $\rho\in M$ a vacuum-like element
if $L(-1)\rho=0$.
\end{defn}

If $\rho$ is a vacuum-like element, we can prove $v_n\rho=0$ for $n\geq 0$ by induction. Because, by the definition of a (weak) $V$-module, there is $m$ such that $v_s\rho=0$ for $s\geq m$. Then $0=L(-1)v_{m}\rho=-mv_{m-1}\rho$ and so
we have $v_n\rho=0$ for $n\geq 0$.
In this case, $V\!\cdot\! \rho=\{v_{-1}\rho\mid v\in V\}$
becomes a $V$-module which is isomorphic to $V$, see Lemma 9.6 in \cite{DLM}.
Here and there, we use a notation $V\!\cdot\! S$
to denote a $V$-submodule generated from a subset $S$
of a $V$-module.

In this section, we will show
the following propositions in a general setting,
which will play important roles
in the proof of Theorem \ref{Main}.

\begin{prop}
Assume that a VOA $V$ of CFT-type has a non-split extension
$$ 0\to V\cong W\to \tilde{V} \xrightarrow{\psi} V \to 0 $$
of $V$-modules with a $V$-homomorphism $\psi$ and $W\cong V$.
Set $\tilde{V}_{[n]}=\{\tilde{v}\in \tilde{V}
\mid \exists m\in \bN \mbox{ such that }(L(0)-n)^m\tilde{v}=0\}$ the genelarized eigenspace of $L(0)$ with eigenvalue $n$.
Let $\tilde{\1}\in \tilde{V}_{[0]}$ be an inverse image of the vacuum element $\1\in V$ by $\psi$ and
set $u=L(-1)\tilde{\1}$.
Then $u\not=0$ and by viewing $u\in V_1$,
$u(0)$ acts on $V_k$ as a nilpotent operator for each $k$.
\end{prop}

\pr
We note $u\in W_1$.
If $u=0$, then $\tilde{\1}$ is a vacuum-like element and
$\tilde{V}$ decomposes into a direct sum of
two $V$-modules $V\cong W$ and $U=\{v_{-1}\tilde{\1}\mid v\in V\}$, which contradicts
the assumption. So, we have $0\not=u\in W\cong V$ 
and we view $u$ as an element of $V$.
Set $P=\{v\in V \mid \exists m\in \bN \mbox{ s.t. }u(0)^mv=0\}$.
Since $u(0)^m(v_kw)=\sum\binom{m}{j}(u(0)^{m-j}v)_k(u(0)^jw)$,
$P$ is a subVA of $V$.
Furthermore, since $\omega(2)u\in V_0=\bC \1$ and $\1(-3)=0$, we have
$$\begin{array}{rl}
u(0)\omega=&u(0)\omega(-1)\1=-[\omega(-1),u(0)]\1\cr
=&
-\{(\omega(0)u)(-1)-(\omega(1)u)(-2)+(\omega(2)u)(-3)+\!\cdots\}\1 \cr
=&\{u(-2)-u(-2)+0\}\1=0.
\end{array}$$
Hence, $P$ is a full subVOA of $V$.
Assume $P\not=V$ and let $\rho$ be a lowest weight homogeneous
element outside of $P$. Then
$P\not\ni u(0)\rho=-\rho(0)u-\sum_{k=1}^{\infty}\frac{(-L(-1))^k}{k!}
\rho(k)u$ by the assumption of $\rho$.
For $k\geq 1$, since $\wt(\rho(k)u)<\wt(\rho)$ and $\omega\in P$,
we have $\rho(k)u\in P$ and
$\sum_{k=1}^{\infty}\frac{(-L(-1))^k}{k!}
\rho(k)u\in P$. Therefore, we have $\rho(0)u\not\in P$.
Viewing it in $V=W\subseteq \tilde{V}$,
$\rho(0)u=\rho(0)\omega(0)\tilde{\1}=\omega(0)\rho(0)\tilde{\1}$.
However, since $\psi(\rho(0)\tilde{\1})=0$ and $\rho(0)\tilde{\1}\in V_{\wt(\rho)-1}$,
we have $\rho(0)\tilde{\1}\in P$. Since $P$ is a full subVOA, we have
$\rho(0)u=L(-1)\rho(0)\tilde{\1}\in P$,
which is a contradiction.
Therefore, $P=V$ and $u(0)$ acts on $V_k$ nilpotently for
every $k\in \bZ$.
\prend

Since $V$ is of CFT-type, $v(0)\tilde{\1}\in W_0=V_0$ for all $v\in V_1$. 
Hence we have 
$$ v(0)u=v(0)L(-1)\tilde{\1}\in L(-1)\bC \1=0$$
for all $v\in V_1$. In particular, $u\in Z(V_1)$, where $Z(V_1)$ denotes the center of 
Lie algebra $V_1$. Therefore, $u, \phi(u)\in Z(V_1)$. 

Since $\langle u,\phi(u)\rangle>0$, we have $\langle u+\phi(u),u+\phi(u)\rangle >0$ or 
$\langle u-\phi(u),u-\phi(u)\rangle>0$.  We can choose $w\in \{\bC(u\pm \phi(u))\}$ such that $\langle w,w\rangle=1$.   
Since $[u(0),\phi(u)(0)]=0$, $w(0)$ also acts on $V_k$ nilpotently.

Let ${\rm VA}(w)$ be a vertex subalgebra generated by $w$. 
Then $\omega^1=\frac{w(-1)w}{2}$ becomes a conformal element of ${\rm VA}(w)$. 
By the relation $[w(m),w(n)]=\delta_{m+n,0}$, we have ${\rm VA}(w)\cong M(1)$ as VOAs, 
where $M(1)$ denotes a VOA of free-boson type with central charge $1$.

\begin{prop}
Under the assumptions in Proposition 1,
if $V$ is a unitary VOA, 
then $V={\rm VA}(w)\otimes \Com({\rm VA}(w),V)$, where 
$\Com({\rm VA}(w), V)=\{v\in V \mid \alpha(n)v=0 \mbox{ for all }\alpha\in {\rm VA}(w) \mbox{ and } n\geq 0\}$.
In particular, $V$ is not holomorphic.
\end{prop}

\pr
Set $Q=\{v\in V\mid w(0)v=0\}$.
As we explained in the proof of Proposition 1, $Q$ is a subVOA and 
$\phi$-invariant, since $\phi(w)\in \pm w$.
Since $V$ has a positive definite form $(,)_V$, $Q^{\perp}$ is also $w(0)$-invariant . 
Since $w(0)$ acts on $(Q^{\perp})_k$ nilpotently for $k\in \bZ$, if $Q^{\perp}\not=0$, then $Q^{\perp}$ contains a nonzero element $v$ such that 
$w(0)v=0$.  
However, since $(,)_V$ is positive definite and $0\not=v\in Q\cap Q^\perp$, 
we have a contradiction. Hence, we have $V=Q$, that is, $w(0)$ acts on $V$ as zero.

Since $w(0)$ acts on $V$ as zero and $V$ has an invariant bilinear form 
$\langle,\rangle$, $V$ is a direct sum of simple $M(1)$-modules $J$ isomorphic 
to $M(1)$, that is, there is $b\in V$ with $\omega^1(1)b=0$ 
such that $J=M(1)b$. 
Since $b\in \Com({\rm VA}(w),V)$,
we have $V\cong {\rm VA}(w)\otimes \Com({\rm VA}(w),V)$. 
We next show that $\omega-\omega^1\in \Com({\rm VA}(w),V)$, where 
$\omega^1=\frac{w(-1)w}{2}$. 
Then $w(0)\omega^1=w(2)\omega^1=0$ and $w(1)\omega^1=w$.
On the other hand, $w(1)\omega=\omega_1w=w$. 
Since $\omega(2)w\in \bC \1$, there is $r\in \bC$ such that $\omega(2)w=r\1$.
Then $r=\langle \omega(2)w,\1\rangle=\langle w, \omega(0)\1\rangle=0$
and so $r=0$, that is, $\omega(2)w=0$.
Therefore, $w(m)(\omega-\omega^1)=0$ for all $m\geq 0$ and 
$\omega-\omega^1\in \Com({\rm VA}(w),V)$.
Since $M(1)$ is not holomorphic, neither is $V$.
\prend

\section{Proof of Theorem \ref{Main}}

We first prove the following result.

\begin{prop}\label{nonsplitextension}
If $V$ is a simple holomorphic irrational VOA of CFT-type, then
there are an indecomposable $V$-module $\tilde{V}$ and its submodule $W$ such that 
$\tilde{V}/W\cong V \cong W$ as $V$-modules. 
\end{prop}

\pr
We start the proof of Proposition \ref{nonsplitextension} by the following case:

\begin{lmm}
Let $V=\coprod_{m=-t}^{\infty}V_m$ be a simple holomorphic VOA of CFT-type and assume that
$A(V)$ is semisimple,
then all $\bN$-gradable $V$-modules $M=\coprod_{n=0}^{\infty}M_{(n)}$ are semisimple. In particular, $V$ is rational.
\end{lmm}

\pr
We note that since $V$ is a holomorphic VOA of CFT-type 
and $A(V)$ is semisimple, $A(V)=\End(V_0)\cong \bC$ 
and $\bar{\omega}=0$, where $\bar{\omega}$ is the image of $\omega$ in $A(V)$. 
Let $M=\coprod_{n=0}^{\infty} M_{(n)}$ be an $\bN$-gradable $V$-module.
Using the notation in \cite{DLM},
$\Omega(M)$ denotes the set of singular elements of a $V$-module $M$, that is,
$\Omega(M)=\{w\in M \mid v_{m}w=0\mbox{ for all }m\geq \wt(v)\mbox{ and }v\in V \}$. 
Since $A(V)$ acts on $\Omega(M)$, $L(0)$ acts on $\Omega(M)$ as zero. 
Furthermore, since $0$ is the lowest eigenvalue of $L(0)$ on a $V$-module,  
$\Omega(M)$ coincides with the generalized eigen subspace of $L(0)$ in $M$ with eigenvalue $0$. 

Since $A(V)\cong \bC$, $\Omega(M)$ is a direct sum of copies of a simple $A(V)$-module $V_0\cong \bC$ 
and furthermore a submodule $V\!\cdot\!\Omega(M)$ of $M$ generated from
$\Omega(M)$ is a direct sum of copies of $V$-module $V$. We will show
$M\subseteq V\!\cdot\! \Omega(M)$. Clearly, $M_{(0)}\subseteq \Omega(M)$ and so $V\!\cdot \!M_{(0)}\subseteq V\!\cdot\! \Omega(M)$.  
If $M_{(1)}\not\subseteq V\!\cdot \!M_{(0)}$, then $M_{(1)}/(V\!\cdot \!M_{(0)}\cap M_{(1)})$ is a top level of $V\!\cdot \!M_{(1)}/(V\!\cdot \!M_{(0)}\cap V\!\cdot \!M_{(1)})$. 
Therefore $A(V)$ acts on $M_{(1)}/(V\!\cdot \!M_{(0)}\cap M_{(1)})$. We hence have $(L(0))M_{(1)}\subseteq  V\!\cdot \! M_{(0)}\cap M_{(1)}$. 
Since the eigenvalues of $L(0)$ is $1$ on $V\!\cdot \!M_{(0)}\cap M_{(1)}$ and $0$ on $M_{(1)}/(V\!\cdot \!M_{(0)}\cap M_{(1)})$, 
we have $M_{(1)}=(\Omega(M)\cap M_{(1)})\oplus (V\!\cdot \!M_{(0)}\cap M_{(1)})$. 
In particular, $M_{(0)}\oplus M_{(1)}\subseteq V\!\cdot\! \Omega(M)$ and the eigenvalues of $L(0)$ on $M_{(1)}$ are $\{0, 1\}$ at most. 
Therefore the eigenvalues of $L(0)$ on $M_{(2)}\cap V\!\cdot\! (M_{(0)}+M_{(1)})$ are $\{1, 2\}$ at most and $0$ on $M_{(2)}/(M_{(2)}\cap V\!\cdot (M_{(0)}+M_{(1)}))$. 
Iterating these steps, we have $M_{(m)}\subseteq V\!\cdot\! \Omega(M)+V\!\cdot\! \sum_{j=0}^{m-1}M_{(j)}$ for all $m$ 
and the eigenvalues of $L(0)$ on $M_{(m)}$ are $\{0,1,...,m\}$ at most. We hence have $M\subseteq V\!\cdot\! \Omega(M)$. 
Since we have already shown that $V\!\cdot\! \Omega(M)$ is a direct sum of simple modules, $M$ is also a direct sum of simple modules. 
\prend

So if $V$ is a simple holomorphic irrational VOA of CFT-type, then $A(V)$ is not semisimple, that is, $J(A(V))\not=0$.

\begin{lmm}
Let $R$ be an associative algebra with an identity $1$ and $J(R)$ a Jacobson radical, that is,
the intersection of all maximal left ideals and assume $R/J(R)\cong
{\rm Mat}_{s,s}(\bC)$.
If $T\subseteq J(R)$ is a left ideal of $R$ such that
$J(R)/T$ has a complement in $R/T$, then $T=J(R)$.
\end{lmm}

\pr
If $R/T=M/T\oplus J(R)/T$, then $R=M+J(R)$ and so there are $m\in M$ and $p\in J(R)$ such that $1=m+p$. In this case, $m=1-p\in M$ is a unit of $R$ and so $M=R$ and $T=J(R)$.
\prend

Viewing $A(V)$ as an $A(V)$-module, set
$M=V\!\cdot\! A(V)=\coprod_{m=0}^{\infty} M_{(m)}$ with $M_{(0)}=A(V)$.
Then we have a non-split exact sequence
$$0\to V\!\cdot\! J(A(V)) \to M \xrightarrow{\psi} V\!\cdot\! (A(V)/J(A(V)))\to 0.$$
We note $A(V)/J(A(V)\cong \bC$ and 
$V\!\cdot\! (A(V)/J(A(V)))\cong V$ as $V$-modules.
Hence $\psi(M_{(k)})=V_k$ for $k\in \bZ$.
Let $\bar{\1}$ be the vacuum element of
$V\!\cdot\! (A(V)/J(A(V)))$ and let $\tilde{\1}$ be an inverse image of $\bar{\1}$ in $M_{(0)}$ and set
$u=L(-1)\tilde{\1}$. Let $W$ be a submodule of 
$V\!\cdot\!J(A(V))$ containing $u$. 
Then, the image in $M/W$ of $\tilde{\1}$ is a vacuum-like element and so
$V\!\cdot\! (A(V)/J(A(V)))$ becomes a direct summand of $M/W$. In particular,
$A(V)/J(A(V))$ is a direct summand of $M_{(0)}/W_{(0)}$.
By the above lemma, we have $W=V\!\cdot\! J(A(V))$. In particular, 
$V\!\cdot\! u=V\!\cdot\!J(A(V))$.
Let $\CF$ be the set of
$V$-submodules of $V\!\cdot\! J(A(V))$ that does not contain $u$.
By Zorn's lemma, $\CF$ has a maximal element, say, $Q$.
Then $V\!\cdot\! J(A(V))/Q=V\!\cdot\! \bar{u}$ and $V\!\cdot\!\bar{u}$ is a simple module, where $\bar{u}$ is the image of $u$ in $M/Q$. Since $M/Q$ is indecomposable module and $M/V\!\cdot\! J(A(V))\cong V$, 
we have a desired indecomposable $V$-module.

This completes the proof of Proposition \ref{nonsplitextension}.
\prend

We now start the proof of Theorem 1. If $V$ is a simple holomorphic irrational VOA of CFT-type,
Then, $V$ satisfies the assumptions in Propositions 1 and 2 by Proposition 3.  
Therefore, $V$ is not holomorphic, which contradicts the assumption.

This completes the proof of Theorem \ref{Main}.

\end{document}